\documentclass[journal]{IEEEtran}
\ifCLASSINFOpdf
\else
\fi
\usepackage[cmex10]{amsmath}

\usepackage{algpseudocode}
\usepackage{algorithm}

\usepackage{amsfonts}
\usepackage{amssymb}

\newtheorem{Definition}{Definition}

\begin{document}
\title{Revisiting the use of Robust Optimization\\for optimal energy offering under price uncertainty}
%
%
%

\author{Fabio~D'Andreagiovanni$^{*}$, 
        Giovanni~Felici, 
        and~Fabrizio~Lacalandra,
\thanks{F.~D'Andreagiovanni is with: the Dept. of Mathematical Optimization, Zuse Institute Berlin (ZIB), Takustr. 7, 14195 Berlin, Germany; the DFG Research Center MATHEON, Technical University Berlin and the Einstein Center for Mathematics (ECMath), Strasse des 17 Juni, 10623 Berlin, Germany, e-mail: d.andreagiovanni@zib.de.}
\thanks{G.~Felici and F.~D'Andreagiovanni are with the Institute for System Analysis and Computer Science ``A. Ruberti'', National Research Council of Italy (IASI-CNR), via dei Taurini 19, 00185 Roma, Italy, e-mail: felici@iasi.cnr.it.}
\thanks{F.~Lacalandra is with Quantek s.r.l., via G. Marconi 3, 40122 Bologna, Italy, e-mail: fabrizio.lacalandra@quantek.it}
\thanks{Manuscript received Month DD, YYYY; revised Month DD, YYYY.}
\\
$^{* \;}$First Author}

%
%

\markboth{Journal of \LaTeX\ Class Files,~Vol.~xx, No.~x, Month~xxxx}%
{Shell \MakeLowercase{\textit{et al.}}: Bare Demo of IEEEtran.cls for Journals}
%



\maketitle

\begin{abstract}
We propose a new Robust Optimization method for the energy offering problem of a price-taker generating company that wants to build offering curves for its generation units, in order to maximize its profit while taking into account the uncertainty of market price.
Our investigations have been motivated by a review of another Robust Optimization method 
proposed in \cite{BaCo11}, which entails the solution of a sequence of robust optimization problems imposing full protection and defined over a sequence of nested subintervals of market prices:
this method presents a number of drawbacks that may sensibly limit its application and computational efficiency in practice and that may expose a company to the risk of presenting offering curves resulting into suboptimal or even infeasible accepted offers. To tackle all such issues, our method provides for solving one single robust counterpart, considering an intermediate level of protection between null and full protection, and for making energy offers at zero price, practically eliminating the risk of non-acceptance. Computational results on instances provided by our industrial partners show that our new method grants a great improvement in profit.
\end{abstract}


%
\IEEEpeerreviewmaketitle

\section{Introduction}
\IEEEPARstart{W}{e} consider the decision problem of a price-taker generating company that operates in a competitive electricity market and wants to select energy offering strategies for its (few) units, in order to maximize the profit while taking into account market price uncertainty. For an introduction to offering in competitive electricity markets, we refer the reader to \cite{DaWe00} and \cite{KwFr12}.

A major challenge in energy offering is represented by tackling uncertainty of market prices: hourly prices forming in the market are not known in advance to companies, which must therefore take offering decisions by making some assumption about future price realizations.
Considering price uncertainty in mathematical optimization models is essential and neglecting it may have very negative consequences and lead to bad offering decisions and losses.
Since the seminal work \cite{Da93}, which analyzes the mathematical consequences of passing from regulated to competitive markets,
several approaches have thus been developed for addressing price uncertainty. We refer the reader to the surveys \cite{LiShQu11} and \cite{TaEtAl15} for a thorough review of related works.

Our main reference for the present work is represented by a \emph{Robust Optimization} method for energy offering introduced by Baringo and Conejo in \cite{BaCo11}. In the remainder of the paper, we will denote such method by \emph{Bar-Con}.
Bar-Con constitutes a recent and relevant reference in literature that builds offering curves by solving a sequence of \emph{(robust) mixed integer linear programming problems}, defined on the basis of a sequence of nested subintervals of market prices. The robust problems follow the well-established $\Gamma$-Robustness model introduced by Bertsimas and Sim  \cite{BeSi04}.
However, Bar-Con presents a number of issues that may sensibly limit its application in practice and expose a company to the risk of building curves that may result in infeasible and suboptimal accepted offers (see Section \ref{sec:BarCon} for our detailed review). Though Bar-Con has been widely cited in other works, we stress that, to the best of our knowledge, so far \emph{no work} available in literature has ever pointed out these issues and tried to overcome them.

\noindent
In this paper, our main original contributions are the following:
\begin{enumerate}
   \item we review the Robust Optimization method proposed by Baringo and Conejo in \cite{BaCo11} highlighting some drawbacks. Specifically, we highlight that:
   \begin{itemize}
   \item
		since Bar-Con imposes full protection against price deviations, the computation of robust optimal solutions does not require to solve robust counterparts including additional 		 
		variables and constraints, as specified by $\Gamma$-Robustness \cite{BeSi04}. We show that the optimal robust solutions corresponding with full protection can be simply obtained by 		
		solving the original problem with modified price coefficients in the objective function.
       This simple observation allows to greatly improve the computational efficiency of the method;
   \item
		offering curves built according to Bar-Con may imply the acceptance of \emph{infeasible} and \emph{suboptimal} offers (e.g., accepted offers may violate ramp limits of the generation 	 
		units, thus being not implementable in practice).
       This is a consequence of building offering curves by merging robust optimal solutions obtained for distinct assumptions about price realizations, as we detail in Subsection
			 \ref{subsec:riskInfeasibilitySuboptimality}.
Moreover, the
method concretely exposes the producer to the
risk that offers are not accepted;
   \end{itemize}
   \item 	we propose a revised method for energy offering under price-uncertainty through Robust Optimization. Our original approach is able to overcome the limits of \emph{Bar-Con} and 	 
					is more computational efficient since it provides for the solution of \emph{only one} $\Gamma$-robust counterpart, based on hourly reference price and price deviations 			 
					derived from historical data.
Our method is more in line with the spirit of $\Gamma$-robustness, reducing the price of robustness by not dealing with extreme conservatism, but addressing intermediate and more rational cases.
A distinctive feature of our approach is to make all selling offers at \emph{zero price}, following the spirit of a pure price-taker that wants to minimize the probability of non-acceptance. Indeed, in current years the growth in zero marginal cost renewable sources has somehow changed the price shape during the hours in day-ahead market. This leads to the risk that in some hours the price falls below the marginal cost of classic fossil fueled power plants, such as Combined Cycle Gas Turbine (CCGT). This risk is however managed by means of the robust approach proposed;
   \item we prove the effectiveness of our revised method on realistic problem instances, considering fossil-fuel power plants and
       referring to price data from the Italian day ahead energy market. 		 Computational tests show that our approach can greatly improve the total profit of the price-taker by considering intermediate level of protection (\emph{in-between} null and full protection).
\end{enumerate}

\smallskip

\noindent
We developed a new method based on Robust Optimization (RO) to exploit the advantages of RO over more traditional optimization under uncertainty approaches, like Stochastic Programming
(see \cite{BeElNe09} and \cite{BeBrCa11} for a thorough discussions).
Computational tractability is in particular a major advantage of RO, deriving from taking into account data uncertainty in a deterministic way, which allows to tackle the so-called ``curse of dimensionality'', typical of
approaches like Stochastic Programming. Such advantages contributed to the wide success of RO in many application fields, including power system optimization, over the last decade (see \cite{TaEtAl15} for an overview).

Besides RO, several other approaches have been proposed over time for tackling price uncertainty in energy offering. Here, we recall just a few main ones, referring the reader to \cite{LiShQu11} for a recent exhaustive survey.
In \cite{CoEtAl04} and \cite{YaSh04}, the price uncertainty is considered in the guise of an objective function that combines the two conflicting objectives of  maximizing profit and minimizing risk through risk tolerance parameters. To reduce the influence of the risk parameters, whose setting may result tricky, \cite{Ja05} proposes to measure risk through conditional value-at-risk. Other works propose to use forecasting techniques to estimate the uncertain prices in the objective function (e.g., \cite{CoEtAl02}). Finally, another distinct body of literature relies on defining sets of price scenarios and following Stochastic Programming approaches
(e.g., \cite{DeEtAl03}, \cite{HrEtAl13}, \cite{LiShLi07}).

\smallskip
\noindent
The paper is organized as follows: in Section \ref{sec:powerOffering}, we introduce the energy offering problem and present an optimization model for it. In Section \ref{sec:BarCon}, we review the Bar-Con method, pointing out its limits. In Section \ref{sec:revisedOffering}, we introduce our revised method based on RO. Finally, in Section \ref{sec:computations}, we present computational results for realistic instances.

\section{The Energy Offering Problem}
\label{sec:powerOffering}
We focus the attention on a generation company that wants to maximize its profit by choosing the output of its generation units in each hour. We assume that the company is a \emph{price-taker}, i.e. it is not able to influence the market price by its generation decisions. An important consequence of this assumption is that the \emph{multi-unit} decision problem can be decomposed into \emph{multiple single-unit} problems, as noted in \cite{BaCo11}.
In a more formal way, we can state the single-unit decision problem as follows:
\begin{Definition}[Energy Offering Problem - EOP]
Given a planning horizon including multiple time periods, the market price in each period and a single generation unit characterized by a set of technical constraints and a generation cost, the \emph{Energy Offering Problem (EOP)} consists in establishing the energy to produce and offer in the market for the unit in each period, in order to maximize the total profit over the time horizon, while satisfying the technical constraints.
\end{Definition}

\medskip
The EOP is an optimization problem that can be modelled as the following \emph{Mixed Integer Programming} problem, which includes continuous and discrete decision variables and which we denote by EOP-MIP
%
\begin{align}
\max
&
\sum_{t \in T} \lambda_{t} p_{t} - c_{t} (p_{t})
&&
\mbox{(EOP-MIP)}
\nonumber
\\
&
p \in P \; ,
&&
\nonumber
\end{align}

\noindent
where:
\begin{itemize}
  \item $T$ is the set of time periods - as in \cite{BaCo11}, we consider hours in a day, so we have $T = \{1,2,\ldots,24\}$;
  \item $p = (p_1, p_2, \ldots, p_{|T|}) \in \mathbb{R}^{|T|}_+$ is a non-negative decision vector, which includes one variable $p_{t} \geq 0$
	    for each period $t \in T$ to represent the energy offered in $t$;
  \item $\lambda_{t}$ is the realized market price in period $t \in T$;
  \item $c_t(p_{t})$ is the cost of generating $p_{t}$ in period $t \in T$, possibly including the so called \emph{Start Up Costs} (SUC);
  \item $P$ is the feasible set of solutions of the problem, containing all vectors $p$ that satisfy the technical constraints of the unit,
				namely minimum and maximum generation output, ramp limits and minimum up and down time.
				We provide a complete description of $P$ in Appendix \ref{sec:appendixRefForm}.
\end{itemize}


\noindent
\textbf{Remark.}
In Appendix \ref{sec:appendixRefForm}, we model the cost function $c_t(p_{t})$ as a quadratic function. In our case, EOP-MIP thus becomes a \emph{Mixed Integer Quadratic Programming (MIQP) problem}.
We note that in literature it is also common to adopt piecewise linear approximations of the quadratic cost function, as done, for example, in \cite{BaCo11}, thus making EOP-MIP a \emph{Mixed Integer Linear Programming} (MILP) problem.
\\
We stress that all the considerations that we make in the next section about defining a robust counterpart for the MIQP version of EOP can be immediately extended to the piecewise linear version considered in \cite{BaCo11}.

\smallskip
\noindent
In the next section, we discuss how to consider the uncertainty of market price in the EOP and we detail a revised and refined robust optimization method that overcomes the limits of the Bar-Con.

\section{A review of the Baringo-Conejo robust offering approach}
\label{sec:BarCon}

In this section, we review the Bar-Con robust optimization-based method and we point out how its computational efficiency can be improved by reducing the size of the robust counterpart. Furthermore, we illustrate how offering curves built according to Bar-Con can lead to suboptimal and even infeasible accepted offers. As first step, we review basic concepts from Robust Optimization and $\Gamma$-Robsutness.

\subsection{A concise introduction to Robust Optimization}

Until this section, we have considered a \emph{deterministic} version of the EOP where we have assumed that the market prices $\lambda_{t}$, $\forall t \in T$ are exactly known when the problem is solved. However, this assumption does not hold in practice: in competitive energy markets, the hourly prices result from the combination of demands and offers whose precise values are not known in advance. The prices $\lambda_{t}$ are thus \emph{uncertain}, i.e. their values are not deterministically known when the EOP is solved, and the price-taker commonly possesses just an estimate of them (e.g., an average value derived from historical market price data, or some other form of prediction).

The presence of uncertain data in the EOP is very tricky:
the price-taker could solve the EOP-MIP referring to the estimates of hourly prices, thus getting solutions that are optimal with respect to such estimates; however, the estimates will be in general (sensibly) different from the actual prices occurring in the market, so that supposedly optimal solutions may reveal to be \emph{heavily sub-optimal}. This fact can be clarified by a simple consideration: if we solve the EOP-MIP using estimates of values $\bar{\lambda}_{t}$ lower than the actual market prices, then we generate less energy than it would be ideal with the real prices, obtaining a lower profit. On the other hand, if we overestimate the actual market prices, we could turn on generation units and generate energy when it is not convenient to do so, facing overproduction and losses.

More in general, the presence of data uncertainty is a source of issues in any optimization problem and neglecting it is well known to have very bad consequences: as we have exemplified above, solutions supposed to be optimal may turn out to be heavily suboptimal and, even worse, solutions supposed to be feasible may instead reveal to violate feasibility constraints, thus resulting in decisions that cannot be implemented in practice. For an exhaustive introduction to the issues arising in optimization under data uncertainty, we refer the reader to the book \cite{BeElNe09}.

In the last decade, RO has known a wide success as a modern methodology to manage data uncertainty in optimization problems, thanks especially to its accessibility and computational tractability. We refer the reader to \cite{BeElNe09} and \cite{BeBrCa11} for an exhaustive introduction to theory and applications of RO.
\\
RO is essentially based on the following principles:
\begin{itemize}
    \item the \emph{actual value} of an uncertain coefficient of the problem is unknown. However, the decision maker is supposed to possess a reasonable reference for the unknown actual value. This reference is called \emph{nominal value} and could be, for example, the expected value of a random variable. The actual value equals the sum of the nominal value and of an unknown deviation;
  \item the decision maker defines an \emph{uncertainty set}. This set characterizes the deviations of coefficients w.r.t. their nominal values, against which the decision maker wants to be protected;
  \item the decision maker derives and solves a \emph{robust counterpart} of its problem: this is a modified version of the original optimization problem that only considers robust feasible solutions, i.e. feasible solutions that are protected against all the deviations of the uncertainty set.
      An \emph{optimal robust solution} is a feasible solution
      granting the best objective value under the worst data deviations;
  \item protection against deviations comes at a cost: the so-called \emph{price of robustness} \cite{BeSi04}. This is a deterioration in the optimal value of a problem that the decision maker must face in order to ensure protection. Uncertainty sets granting higher protection entail in general a higher price of robustness.
\end{itemize}

\noindent
Over the years, many models have been proposed to represent uncertainty sets in RO (see \cite{BeBrCa11} for an overview). Among them, one of the most successful has been the $\Gamma$-Robustness model (\emph{$\Gamma$-Rob}) proposed by Bertsimas and Sim \cite{BeSi04}, which has been also adopted in Bar-Con. We  thoroughly review the adaption of $\Gamma$-Rob in Bar-Con in the next subsection.
Here, we just recall the essence of this classical and widely used robust optimization model. Given an uncertain constraint of the problem including $n$ uncertain coefficients of which we know the nominal value and the maximum deviation w.r.t. the nominal value, $\Gamma$-Rob provides for an elegant theory to define a robust counterpart for an uncertainty set imposing full protection against at most $0 \leq \Gamma \leq n$ coefficients deviating from their nominal value. The  $\Gamma$-Rob counterpart has the desirable property of entailing a ``contained'' increase in the number of constraints and variables of the problems (more formally, the counterpart is \emph{compact}) and of maintaining the features of the original problem (e.g., an uncertain MILP problem has a MILP $\Gamma$-Rob counterpart, whereas an uncertain MIQP has a MIQP $\Gamma$-Rob counterpart).
The parameter $\Gamma$ controls the robustness: for $\Gamma=0$ no protection is imposed and no price of robustness occurs; for increasing $\Gamma$, the protection and the price of robustness increase, until for $\Gamma=n$ the highest protection is reached.  This informal description is detailed in the next paragraph with reference to the EOP.

\subsection{The Robust Optimization approach by Baringo and Conejo}
\label{SubSect:BaringoConejoApproach}

In this section, we provide a description of the main features of the Bar-Con, which are relevant to highlight its limits and derive our improved method. For a full description of the method, we refer the reader to \cite{BaCo11}. The Bar-Con assumes that in every hour $t \in T$ the actual market price $\lambda_{t}$ lies in a known range, i.e. $\lambda_{t} \in [\lambda_{t}^{\min}, \lambda_{t}^{\max}]$, and a-priori defines an elementary market price shortfall equal to $\delta \hspace{0.05cm} (\lambda_{t}^{\max} - \lambda_{t}^{\min})$ with $0 \leq \delta \leq 1$. Then, it proposes to compute the hourly offering curves by solving a sequence of $1,\ldots, K$ robust counterparts of the original problem and in the generic k-th robust counterpart:
\begin{enumerate}
  \item the market price in each hour is supposed to lie in the interval
    $$[\lambda_{t}^{\max} - d_{t}^k, \hspace{0.2cm} \lambda_{t}^{\max}]
    $$
    where $d_{t}^k = G^k  (\lambda_{t}^{\max} - \lambda_{t}^{\min}) =  (k - 1) [\delta (\lambda_{t}^{\max} - \lambda_{t}^{\min})]$ is the \emph{maximum price deviation} in $t$ at iteration $k$. Such deviation
causes the overall deviation interval $[\lambda_{t}^{\max} - d_{t}^k, \hspace{0.1cm} \lambda_{t}^{\max}]$ to become wider as the iteration $k$ increases. We remark that \cite{BaCo11} just refers to the coefficient $G^k$ and does not use its equivalent form $(k - 1) \delta$ that we adopt here to highlight the role of the iteration index $k$ in the product;
  \item data uncertainty in the problem is modelled through $\Gamma$-Robustness and
  the k-th robust counterpart has the following form:
\begin{eqnarray}
\max
&&
\sum_{t \in T} \left[ \lambda_{t}^{\max} p_{t} - c_{t} (p_{t}) \right] - \Gamma z - \sum_{t \in T} q_{t}
\nonumber
\\
&&
z + q_{t} \geq d_{t}^k y_{t}
\hspace{0.70cm}
t \in T
\nonumber
\\
&&
z \geq 0
\nonumber
\\
&&
q_{t} \geq 0
\hspace{1.80cm}
t \in T
\nonumber
\\
&&
p_{t} \leq y_{t}
\hspace{1.65cm}
t \in T
\label{linkingConstraints-BarCon}
\\
&&
p \in P
\nonumber
\end{eqnarray}
with $\Gamma$ set equal to the number of periods $|T|$.
\item offering curves for each period are built merging the $K$ optimal robust solutions: the optimal robust solution of the k-th counterpart specifies the generation output $p_t^k$ to be offered at price $\lambda_{t}^{\max} - d_t^k$. For each period $t \in T$, the merging of prices and generation volumes for all the iterations $k \in K$ defines a non-decreasing offering curve.
\end{enumerate}

\noindent
In the next paragraphs, we review the Bar-Con approach highlighting its limits, which may severely reduce the advantages of using it in practice. We first show by a few simple observations that the robust counterpart can be greatly simplified. Specifically, we show that: i) the variables $y_t$, $t \in T$ and the constraints $\eqref{linkingConstraints-BarCon}$ are not necessary in the robust counterpart; ii) the robust counterpart can be just defined as a modified version of the original problem EOP-MIP and the robust auxiliary variables $q_t$, $z$ are actually not necessary.
Secondly, we show that offering curves defined according to the Bar-Con method may lead to infeasible and suboptimal accepted offers.
Last but not least, we note that this approach provides for building stepwise offering curves with a high number of steps and would thus result not adequate for many markets, where the number of steps that define a curve is bounded (for example, in the day ahead Italian Energy Market the limit is 4, see \cite{GME_website}).

\subsection{Reducing the size of the robust counterpart}
\label{subsec:reducingSize}

As first step of our review of Bar-Con, we show that the robust counterpart does not need the auxiliary variables $y_t$ and the additional constraints linking $y_t$ and $p_t$. We show this referring to the canonical passages used in \cite{BeSi04} to get the $\Gamma$-robust counterpart of an optimization problem.
We first rewrite the robust version of EOP-MIP at iteration $k$ as:
\begin{align}
\max
&
\sum_{t \in T} \left[ \lambda_{t}^{\max} p_{t} - c_{t} (p_{t}) \right] - DEV_k (\Gamma,p)
&&
\nonumber
\\
&
p \in P
&&
\nonumber
\end{align}

\noindent
where the additional term $- DEV_k (\Gamma,p)$ represent the worst reduction that the objective function may experience. Such term is defined for a power vector $p$, when at most $\Gamma$ hourly prices deviate from their reference value $\lambda_t$ at iteration $k$. The value $DEV_k (\Gamma,p)$ corresponds to the optimal value of the following binary program:
\begin{align}
DEV_k (\Gamma,p) = \max
&
\sum_{t \in T} (d_{t}^k p_{t}) \hspace{0.1cm} w_{t}
&&
\nonumber
\\
&
\sum_{t \in T } w_{t} \leq \Gamma
&&
\nonumber
\\
&
w_{t} \in \{0,1\}
&&
t \in T \; .
\nonumber
\end{align}

In this problem, 1) a binary variable $w_{t}$ is equal to 1 if the price in hour $t$ deviates from its nominal value and experiences the worst deviation $d_{t}^k p_{t}$, whereas it is equal to 0 otherwise; 2) the single constraint imposes an upper bound $0 \leq \Gamma \leq |T|$ on the number of hours whose prices may deviate; 3) the objective function maximizes the deviation from the nominal value for a given output $p_t$, $t \in T$.

The robust version of EOP-MIP with $DEV_k (\Gamma,p)$ presents an inner maximization problem which includes the products $p_{t} \hspace{0.1cm} w_{t}$ of decision variables. However, this is not a real issue since, as proved in \cite{BeSi04}, such non-linearities can be linearized leading to a compact robust counterpart.

First, we note that for a fixed vector $p_t$, the value $DEV_k (\Gamma,p)$  is equal to the optimal value of its \emph{linear relaxation}:
\begin{align}
DEV (\Gamma,p) = \max
&
\sum_{t \in T} (d_{t}^k p_{t}) \hspace{0.1cm} w_{t}
&&
\nonumber
\mbox{(DEV-primal)}
\\
&
\sum_{t \in T } w_{t} \leq \Gamma
&&
\nonumber
\\
&
0 \leq w_{t} \leq 1
&&
t \in T \; .
\nonumber
\end{align}

\noindent
We can then define the \emph{dual problem} of the previous linear program, i.e.:
\begin{align}
\min
\hspace{0.2cm}
&
\Gamma z + \sum_{t \in T} q_{t}
&&
\mbox{(DEV-dual)}
\nonumber
\\
&
z + q_{t} \geq d_{t}^k p_{t}
&&
t \in T
\nonumber
\\
&
z \geq 0
&&
\nonumber
\\
&
q_{t} \geq 0
&&
t \in T \; .
\nonumber
\end{align}

\noindent
Since DEV-primal is feasible and bounded, by strong duality also its dual problem DEV-dual is feasible and bounded and their optimal values coincide. We can then substitute the term $DEV_k (\Gamma,p)$ in our original robust counterpart with the dual problem as follows:
\begin{align}
\max
&
\sum_{t \in T} \left[ \lambda_{t}^{\max} p_{t} - c_{t} (p_{t}) \right] - \Gamma z - \sum_{t \in T} q_{t}
&&
\mbox{(Rob-EOP)}
\nonumber
\\
&
z + q_{t} \geq d_{t}^k p_{t}
&&
t \in T
\label{Rob-EOP_dualCnstr}
\\
&
z \geq 0
&&
\label{Rob-EOP_dualVar1}
\\
&
q_{t} \geq 0
&&
t \in T
\label{Rob-EOP_dualVar2}
\\
&
p \in P \; .
&&
\nonumber
\end{align}

In contrast to the first robust counterpart including the term $DEV (\Gamma,p)$, this robust counterpart does not include the inner maximization problem and the product of variables. Additionally, the increase in size due to the presence of the additional dualization variables and constraints (\ref{Rob-EOP_dualCnstr}-\ref{Rob-EOP_dualVar2}) is ``contained'': the formulation is indeed \emph{compact}, i.e. its size is polynomial in the size of the input.

We emphasize that the adopted dualization procedure is nothing but the standard one associated with $\Gamma$-Rob and proposed in \cite{BeSi04}. As evident, it does not need the presence of the auxiliary variables $y_{t}$ and of the constraints $p_t \leq y_t$. The size of the robust counterpart is thus reduced by eliminating $|T|$ constraints and $|T|$ decision variables. The elimination of these variables ancd constraints represents  an advantage for solving the problem for any commercial optimization solver.

To conclude, we note that when EOP-MIP includes a quadratic cost function and is thus an MIQP, Rob-EOP is an MIQP, whereas when EOP-MIP includes a piecewise linear cost function and is thus an MILP, Rob-EOP is an MILP.

\subsection{Full robustness as a cost-modified nominal problem}
\label{subsec:fullRobustness}

The linear and compact robust counterpart Rob-EOP depends upon the parameter $\Gamma$, which indicates the number of price coefficients against the deviation of which the decision maker wants to be protected. In the case of EOP, $\Gamma$ corresponds with daily hours and can range from 0 to 24.
When $\Gamma \in \{ 1, \ldots, 23\}$, a robust optimal solution can be obtained by solving Rob-EOP.
When instead $\Gamma = 0$ or $\Gamma = 24$, we are considering extreme cases: if $\Gamma = 0$, then no hourly price is assumed to deviate, there is no protection, and we get the original deterministic problem; in contrast, if $\Gamma = 24$, then all the coefficients are assumed to deviate, we get full protection, and we do not need to solve Rob-EOP. Here, we stress that it is indeed just sufficient to consider the following problem:
\begin{align}
\max
&
\sum_{t \in T} \left(\lambda_{t}^{\max} - d_{t}^k \right) - c_{t} (p_{t})
&&
\mbox{(Worst-EOP)}
\nonumber
\\
&
p \in P \; ,
&&
\nonumber
\end{align}

\noindent
which is a version of EOP that considers the worst price $\lambda_{t}^{\max} - d_{t}^k$ that  may occur at iteration $k$ in every hour. Since Bar-Con imposes $\Gamma = 24$ in each of the $k$ robust counterparts (see the algorithm flowchart of Figure 2 in \cite{BaCo11}), then the use of the dualization procedure of $\Gamma$-Rob can be completely avoided. The worst-case robust counterpart Worst-EOP thus eliminates the need for the dualization variables and constraints (\ref{Rob-EOP_dualCnstr}-\ref{Rob-EOP_dualVar2}), reducing the size of the counterpart by $|T|+1$ variables and $|T|$ constraints. This has the effect of improving the efficiency of the branch-and-cut-based solution process implemented by commercial optimization solvers.

\subsection{Risk of building offering curves leading to infeasible and suboptimal offers}
\label{subsec:riskInfeasibilitySuboptimality}

The offering curves built according to Bar-Con are obtained merging $K > 0$ robust optimal solutions obtained for $K$ distinct configurations of the prices over the 24-hour horizon, as we explained in Section \ref{SubSect:BaringoConejoApproach}.

Merging optimal solutions concretely exposes the generating company to the risk of building offering curves that may result into \emph{infeasible} and
\emph{suboptimal} accepted offers. We illustrate this fact by a simple realistic example, defined with our industrial partners.

\smallskip
\noindent
\textbf{Example 1}. We consider the status and the output of a generation unit in a time interval made up of 3 consecutive hours that we conventionally denote by indices $t = 1,2,3$, assuming that in an initial reference hour $t = 0$ the unit is off and the  output is null (i.e., $p_0 = 0$). The generation unit has a minimum power output of 160 MW and a capacity of 440 MW with a startup ramp limit of 160 MW/h and ramp-up limit of 55 MW/h. The production cost is expressed by a quadratic function of the power output that ranges from a value of 8768 EUR for the minimum output
to 25848 EUR for the maximum output.
Referring to real prices of the Italian energy market, we assume that the maximum price observed in the 3 hours are $\lambda^{\max}_1 = 54$, $\lambda^{\max}_2 = 55$, $\lambda^{\max}_3 = 61$ (in EUR/MWh). These prices define a market price configuration that we denote by C1. Following the approach Bar-Con, from C1 we can define two other price configurations C2 and C3, by reducing the prices of C1 by 1 and 2 EUR respectively (see Table \ref{table:prices})

We then solve the robust counterpart Worst-EOP corresponding with full protection $\Gamma = 24$, obtaining the robust optimal outputs shown in Table \ref{table:output}. While for the maximum price configuration C1, the optimal choice is to ramp up to the maximum in each period, when the prices decrease it is instead optimal to keep turned off the unit for an increasing number of hours: for C2, the unit is kept off in hour 1, whereas for C3 it is kept off in hours 1 and 2.
Following Bar-Con, the robust optimal solution leads to the definition of the hourly offering curves shown in Table \ref{table:offeringCurves}, where each curve is actually a 3-step function with each step associated with a pair (price, offered output).
\begin{table}
\renewcommand{\arraystretch}{1.3}
\caption{Hourly price configurations}
\label{table:prices}
\small
\begin{center}
\begin{tabular}{c c c c}
\hline
Configuration  & $\lambda_1$ 	& $\lambda_2$ & $\lambda_3$
\\ [2pt]
\hline
C1 	&  	54	&  	55	& 	61\\
C2 	&  	53	&  	54	& 	60\\
C3 	&  	52	&  	53	& 	59\\
\hline
\end{tabular}
\end{center}
\end{table}

\begin{table}
\renewcommand{\arraystretch}{1.3}
\caption{Robust optimal output and profit for each price configuration}
\label{table:output}
\small
\begin{center}
\begin{tabular}{c c c c c}
\hline
Configuration  & $p_1$ 	& $p_2$ & $p_3$ & Total profit
\\ [2pt]
\hline
C1 	&  	160	&  	215	& 	270 & 1498\\
C2 	&  	0	&  	160	& 	215 & 1020\\
C3 	&  	0	&  	0	& 	160 & 672\\
\hline
\end{tabular}
\end{center}
\end{table}

\begin{table}
\renewcommand{\arraystretch}{1.3}
\caption{Hourly offering curves}
\label{table:offeringCurves}
\small
\begin{center}
\begin{tabular}{c c c c}
\hline
t  & Step 1 & Step 2 & Step 3
\\ [2pt]
\hline
1 	&  	(0, 52)	&  	(0, 53)	&  (160, 54)\\
2 	&  	(0, 53)	&  	(160, 54)	& (215, 55)\\
3 	&  	(160, 59)	&  	(215, 60)	& (270, 61)\\
\hline
\end{tabular}
\end{center}
\end{table}

Since the price-taker specifies the minimum selling price for a given quantity, it may happen that offering curves will not be (entirely) accepted, thus leading to generation plans that violate the ramp constraints over consecutive hours. Referring to Table \ref{table:offeringCurves}, this happens, for example, if the actual price occurring in the three hours are $52$, $53$ and $61$, respectively: in this case, in hour t=1 and t=2 the market operator ``accepts'' the 0 quantity offers of the price-taker, while in hour t=3 it accepts the offer equal to 270 MWh. Since the unit of the price taker cannot ramp from 0 to 270 in 1 hour, the offering is evidently infeasible in practice.

Another issue with Bar-Con is represented by the risk of suboptimality of offering. This issue again derives from the fact that Bar-Con provides for specifying a minimum price that is requested to sell a specific output. To illustrate the risk of suboptimality, we suppose that the actual prices occurring in the three hours are $54$, $53$, and $59$: the accepted offers of the price-taker are in this case $(160, 0, 160)$. This solution is feasible w.r.t. the ramp constraints. However, it is associated with a total profit of 544 EUR, that is sensibly lower than the optimal profit of 672 EUR that would be obtained for the optimal production equal to (0, 0, 160). In this case, the suboptimality derives from turning on the unit
in $t=1$, action that would not be taken in an optimal solution. 

\section{A revised $\Gamma$-robustness approach for energy offering}
\label{sec:revisedOffering}

In this section, we illustrate our new method for energy offering, which tackles the drawbacks of Bar-Con that we have highlighted previously. Our method presents three main advantages:
\begin{enumerate}
  \item it properly manages the risk of defining offers that are not accepted in the market;
  \item it tackles the issue of defining offering curves that result into infeasible and suboptimal accepted offers;
  \item it is in line with the spirit of Robust Optimization and $\Gamma$-Robustness, according to which full protection (i.e., $\Gamma = 24$ as in Bar-Con) should be avoided (see \cite{BeElNe09,BeBrCa11}); computationally,  our method indicates that intermediate protection (i.e., $0 < \Gamma < 24$) guarantees higher profits.
\end{enumerate}

\noindent
To achieve this, we propose three strategies unified in a new robust optimization method: 
\begin{enumerate}
  \item all our offers are made at \emph{zero price} (i.e., we do not require a minimum positive price for selling energy), so that the offers will be always less or equal than the market clearing price (under the assumption that we operate in an energy market where negative prices are not allowed and the lowest feasible price is 0 -  such as in the Italian Energy Market \cite{GME_website}). The risk of non acceptance bidding at zero price thus becomes negligible;
  \item we compute a single $\Gamma$-robust optimal solution referring to one single price configuration derived from historical data. So we avoid the infeasibility and suboptimality issues associated with merging distinct robust optimal solutions obtained for multiple price configurations;
  \item we define the single robust optimal solution by solving the robust counterpart with $\Gamma$ such that $0 < \Gamma < 24$, so to avoid the excessive conservatism associated with full protection.
\end{enumerate}

The robust counterpart that we adopt differs from that of Bar-Con concerning the nominal prices and the worst price deviations used in it, which we define by using realizations of the market price occurred in previous days. To this end, we suppose that for each $t \in T$ we have a number $I>0$ of observed past realizations of the price $\lambda_t^{i}$ with $i = 1,\ldots,I$. Without loss of generality, we assume that these values are sorted in non-decreasing order, i.e. $\lambda_t^{i} \leq \lambda_t^{i+1}$ for $i = 1,\ldots,I-1$. The assumption of possessing past price realizations is fully realistic, since hourly prices are typically promptly disclosed by the market operator.
Using these past realizations, we define:
\begin{itemize}
  \item the \emph{nominal price} $\lambda_t^{\text{NOM}}$ for each $t \in T$ as the average value of past observations in $t$, namely:
      $$
      \lambda_t^{\text{NOM}} \hspace{0.1cm} = \hspace{0.1cm}
      \lambda_t^{\text{AVG}} \hspace{0.1cm} = \hspace{0.1cm}
      \frac{1}{I} \hspace{0.1cm} \sum_{i = 1}^{I} \lambda_t^{i}
      $$
  \item the \emph{worst deviation} $d_t$ used in each $t \in T$ as the difference
      $$
      d_t \hspace{0.1cm} = \lambda_t^{\text{AVG}} - \hspace{0.1cm} \lambda_t^{J+1}
      $$
      where $\lambda_t^{J+1}$ is the smallest price observation obtained after excluding the $0 \leq J < I$ smallest  observations (since we have assumed that observations are sorted in non-decreasing order, this means that we are excluding observations $\lambda_t^{i}$ with $i = 1, \ldots, J$ and the worst ``significative'' observation is then $\lambda_t^{J+1}$).
\end{itemize}
We note that any other predictive procedure to obtain nominal prices and deviation can be considered in our framework as a plug-in strategy.
For instance, one could use a Machine Learning approach or a simulation model.
We note that this choice of the nominal value has been inspired by a practice that we have observed among power systems professionals: offering  decisions are taken considering the average price for each hour in a previous time horizon. Also, the choice of excluding a number of the worst price realizations, commonly referred to as {\it trimming}, is widely use in practice for the management of outliers in distributions, and has been suggested and validated by our industrial partners.

\noindent
Our reference robust counterpart then becomes:
\begin{align}
\max
&
\sum_{t \in T} \left[ \lambda_t^{\text{AVG}} p_{t} - c_{t} (p_{t}) \right] - \Gamma z - \sum_{t \in T} q_{t}
&&
\text{(MOD-Rob-EOP)}
\nonumber
\\
&
z + q_{t} \geq (\lambda_t^{\text{AVG}} - \lambda_t^{J+1}) \hspace{0.1cm} p_{t}
&&
t \in T
\nonumber
\\
&
z \geq 0
&&
\nonumber
\\
&
q_{t} \geq 0
&&
t \in T
\nonumber
\\
&
p \in P \; ,
&&
\nonumber
\end{align}

\noindent
which is the canonical $\Gamma$-Robustness counterpart that we have derived in Section
\ref{SubSect:BaringoConejoApproach}, except for the setting of the nominal price and the worst price deviation.

Our new method provides for solving a single time the robust counterpart MOD-Rob-EOP, by making one single assumption about the price configuration. The obtained robust optimal solution is then used for establishing our offering in the market as we detail in Algorithm \ref{ALGrobustOffering}.
\begin{algorithm}
\caption{{\sc $\Gamma$-Offering}}
\label{ALGrobustOffering}
\begin{algorithmic}[1]
    \Require
        Past hourly price realizations $\lambda_t^{i}$,
        $\forall i = 1,\ldots,I$, $\forall t \in T$;
        number $J$ of smallest price realizations to exclude;
        protection parameter $\Gamma:$ $0 \leq \Gamma \leq 24$
    \Ensure Hourly power quantities $p_t^{\star}$ to be offered at zero price
    \State Compute the average hourly prices $\lambda_t^{\text{AVG}}$ $\forall t \in T$, using observations $\lambda_t^{i}$
    %
    \State Solve the robust counterpart MOD-Rob-EOP for $\Gamma$ using:
    \begin{itemize}
      \item $\lambda_t^{\text{AVG}}$ as nominal value $\forall t \in T$
      \item the difference $\lambda_t^{\text{AVG}} - \lambda_t^{J+1}$ as worst  deviation $\forall t \in T$
    \end{itemize}
    \State Let $p^{\star}$ be a robust optimal output vector for MOD-Rob-EOP, then $p_t^{\star}$ is the quantity to be offered in the market at zero price in each period $t \in T$
\end{algorithmic}
\end{algorithm}

A critical difference of our new method formalized through the algorithm $\Gamma$-Offering is that we solve one single $\Gamma$-robust counterpart defined for one single assumption on prices. This is in contrast to Bar-Con where a large number of $\Gamma$-robust counterparts is provided to be solved (in the experimental section of \cite{BaCo11}, the number of solved counterparts is $k = 100$).

As we comment in the following section, the solution of the robust counterpart can be made almost instantaneously, thus making our new method fully suitable for being used in the daily offering decision process of a price-taker company.

\section{Computational results}
\label{sec:computations}
We assessed the performance of our new method for robust energy offering by considering a set of 45 instances provided by our industrial partners.
The 45 instances are based on 15 generation units of different size and features located in three distinct price zones of Italy. The capacity of the units ranges from about 100 to 1250 MW.
The price coefficients used in the robust counterparts are derived from real prices occurred in the Italian day-ahead Market (\emph{Mercato del Giorno Prima (MGP)} \cite{GME_website}) from January 1st to December 31st 2014.

Our computational tests have
one major objective: showing that the solution of a robust counterpart imposing a (small) intermediate level of protection can grant a very big improvement in profit with respect to the robust counterparts imposing null and full protection (i.e., the strategy that we have observed among professionals and the strategy protecting against the worst price deviation in each hour, resembling the protection imposed by Bar-Con).
We note that in our experiments we do not use the method Bar-Con as a benchmark, because of its limits that we have highlighted in Section \ref{SubSect:BaringoConejoApproach}
and that concretely expose a company to the risk of suboptimal and infeasible accepted offers.

We performed all the experiments on a 2.70 GHz machine with 8 GB. The code was written in the C/C++ programming language and the optimization problems were solved by IBM ILOG CPLEX 12.5 interfaced with the code through Concert Technology. We stress that solving MOD-Rob-EOP to optimality does not constitute at all a challenge for CPLEX: for all instances and for all the values of $\Gamma$, CPLEX is able to return an optimal solution within 1 second, making the method fully suitable for daily energy offering.

For each generation unit and each hour $t \in T$, the average price $\lambda_t^{\text{AVG}}$ and the worst deviation $\lambda_t^{\text{AVG}} - \lambda_t^{J+1}$ are defined considering 20 realizations of prices occurring in hour $t$ and in the price zone of the unit, in a time window of 4 consecutive weeks excluding Saturdays and Sundays. We consider three percentage of exclusions of the smallest price realizations: 0, 10 and 20\%. This leads to respectively exclude the smallest $J = \{0, 2, 4\}$ realizations.
Increasing $J$ corresponds with considering a less risk-adverse price-taker that neglects a higher number of extreme deviations.
Such data are used as input for defining the robust counterpart for the 4-week time window, whose optimal robust solution is evaluated in the week following the 4 weeks of the window (excluding again Saturday and Sunday), as we explain in the next paragraph.
The overall (4+1)-week time window that we consider
is shifted by 1 week through the entire 2014,
from the first week of January till the last week of November. This gives raise to 24 time windows that we denote by $W = \{1,2,\ldots,24\}$.

For each of the 15 units, we solve one robust counterpart Mod-Rob-EOP for each value of $J \in \{0, 2, 4\}$, value of $\Gamma \in \{0,1,\ldots,24\}$ and time window $W \in \{1,2,\ldots,24\}$, obtaining a robust optimal output $p^{\star}(J,\Gamma,W)$ (note that here we highlight the dependency of $p^{\star}$ upon the values of $J$ and $\Gamma$ and the window W). The vector $p^{\star}(J,\Gamma,W)$ specifies our \emph{zero-price} energy offering for the unit in the test week of window $W$. We evaluate the performance of each $p^{\star}(J,\Gamma,W)$ in the test week of $W$ as follows:
in each hour $t \in T$ of each day, $p^{\star}_t(J,\Gamma,W)$ is offered giving a total profit over the week that we denote by $\pi(J,\Gamma,W)$. We then sum such profit over all the 24 time windows in 2014 getting a total yearly profit equal to $\pi(J,\Gamma) = \sum_W \pi(J,\Gamma,W)$.

The complete results of our computational tests are presented in Table \ref{table:results}, where we denote by $\Gamma_{BEST}$ the value of $\Gamma$ that grants the highest profit for a given generation unit and value of $J$ for a window. In the table, \emph{ID} identifies the generation unit (we note that the size of the units increases as the ID increases) and \emph{\%Ex} is the percentage of price realizations that are excluded from the input data of the instance. The following four columns
$\Delta \pi$, $\Delta \pi \%$ for $\Gamma_{BEST}$ w.r.t. $\Gamma_0$
and $\Delta \pi$, $\Delta \pi \%$ for $\Gamma_{BEST}$ w.r.t. $\Gamma_{24}$
report the absolute (in euros) and the percentage increase that the robust optimal solution obtained for $\Gamma_{BEST}$ grants with respect to $\Gamma = 0$ and $\Gamma = 24$, respectively.

We stress that for all instances, $\Gamma_{BEST}$  assumes a small value between 1 and 4 and \emph{never} coincides with $0$ (i.e., no protection) and $24$ (i.e., full protection).
This has a profound impact when comparing the profit of $\Gamma_{BEST}$ to that of $\Gamma_{24}$ (reported in the 5th and 6th column of Table \ref{table:results}): setting $\Gamma = 24$ and thus imposing protection against the worst price deviation  in each hour leads to extremely conservative solutions that dramatically decrease the profit obtainable by $\Gamma_{BEST}$. The reduction in profit happens in all 90 but 3 cases and ranges from about $100\,000$ to even more than 10 millions of euros. The over-conservatism is particularly evident when comparing the best robust profit to that obtained for full protection (column $\Gamma_{BEST}$ w.r.t. $\Gamma_{24}$): in the case of instances corresponding to units of higher capacity (i.e, U11 - U15), $\Gamma_{BEST}$ grants an increase in profit that in most cases is of several millions, reaching a peak of 12 millions for instance U12 and $J = 0$. This clearly indicates through experiments a cardinal fact associated with the theory of $\Gamma$-Robustness: imposing full protection and thus setting $\Gamma$ to the highest possible value should be absolutely avoided since it generally entails protection against extreme unlikely cases at the cost of an extremely high price of robustness.
In contrast, as we showed trough the computational results, small values of $\Gamma$ are able to guarantee protection against price uncertainty, while ensuring extremely satisfying level of profits.


%
\begin{table}
\renewcommand{\arraystretch}{1.3}
\caption{Computational results}
\label{table:results}
\small
\begin{center}
\begin{tabular}{cc|rr|rr}
\hline
 &  	& \multicolumn{2}{c|}{\boldmath$\Gamma_{BEST} \mbox{ \textbf{w.r.t.} } \Gamma_0$}
& \multicolumn{2}{c}{\boldmath$\Gamma_{BEST} \mbox{ \textbf{w.r.t.} } \Gamma_{24}$}
\\
\textbf{ID} & $\textbf{\%Ex}$ 	& $\Delta \pi (EUR)$ 	& $\Delta \pi \%$
& $\Delta \pi (EUR)$ &  $\Delta \pi \%$
\\ [2pt]
\hline
 	&  	0   &  + 40,400     &  + 5.8 &  + 213,731   & + 40.5   \\
U1 	&  	10  &  + 44,350     &  + 6.3 &  + 183,161   & + 32.5   \\
 	&  	20  &  + 41,922     &  + 6.0 &  + 152,608   & + 25.8   \\
\hline
 	&  	0     & + 23,395   &   + 5.0 &  + 333,543   & + 212.1\\
U2 	&  	10    & + 46,064   &   + 9.9 &  + 234,371   & + 84.0\\
 	&  	20    & + 42,072   &   + 9.0 &  + 218,607   & + 75.2\\
\hline
  	&  	0     & - 1,383   &   - 0.1 &  + 1,984,511   & + 47.6\\
U3 	&  	10    & + 88,980  &   + 1.4 &  + 1,031,466   & + 19.8\\
 	&  	20    & + 105,253 &   + 1.7 &  + 627,146    &  + 11.1\\
\hline
  	&  	0    & + 43,246  &    + 6.3 &  + 255,124    & + 53.5\\
U4 	&  	10   & + 57,387  &    + 8.3 &  + 181,357    & + 32.1\\
 	&  	20   & + 51,615  &    + 7.5 &  + 148,635    & + 25.1\\
\hline
  	&  	0    & + 15,454  &    + 3.6 &   + 340,568    & + 319.0\\
U5 	&  	10   & + 45,328  &    + 10.5 &  + 240,506    & + 101.6\\
 	&  	20   & + 45,331  &    + 10.5 &  + 199,407    & + 71.8\\
\hline
  	&  	0    & + 14,274  &    + 5.3 &   + 2,030,185    & + 44.9\\
U6 	&  	10   & + 91,766  &    + 10.6 &  + 1,117,143    & + 20.2\\
 	&  	20   & + 152,707  &   + 11.8 &  + 675,172    & + 11.2\\
\hline
  	&  	0    & + 307,691  &    + 5.7 &   + 1,312,795   & + 30.1\\
U7 	&  	10   & + 268,509  &    + 5.0 &   + 909,989    & + 19.3\\
 	&  	20   & + 195,207  &    + 3.6 &   + 792,081    & + 16.6\\
\hline
  	&  	0    & + 465,184  &    + 12.6 &   + 1,295,788   & + 45.4\\
U8 	&  	10   & + 492,569  &    + 13.4 &   + 1,052,153   & + 33.7\\
 	&  	20   & + 387,576  &    + 10.5 &   + 888,434    & + 27.9\\
\hline
  	&  	0    & - 179,097  &    - 0.5 &   + 8,606,256   & + 35.5\\
U9 	&  	10   & + 249,549  &    + 0.8 &   + 4,586,552   & + 16.0\\
 	&  	20   & + 253,872  &    + 0.8 &   + 3,008,977   & + 9.9\\
\hline
  	&  	0    & + 579,613  &    + 11.7 &   + 1,711,145  & + 44.8\\
U10 &  	10   & + 662,118  &    + 13.4 &   + 870,430   & + 18.3\\
 	&  	20   & + 502,539  &    + 10.1 &   + 594,417   & + 12.2\\
\hline
  	&  	0    & + 612,087  &    + 19.6 &   + 2,471,071  & + 196.0\\
U11 &  	10   & + 465,015  &    + 14.9 &   + 1,270,815   & + 54.9\\
 	&  	20   & + 534,280  &    + 17.1 &   + 788,373   & + 27.5\\
\hline
  	&  	0    & + 23,936  &    + 0.1 &   + 12,397,480  & + 49.4\\
U12 &  	10   & + 409,220  &    + 1.1 &  + 4,751,122   & + 14.3\\
 	&  	20   & + 438,452  &    + 1.2 &  + 3,506,378   & + 10.2\\
\hline
  	&  	0    & + 111,222  &    + 0.3 &  + 8,244,776  & + 33.6\\
U13 &  	10   & + 421,917  &    + 1.3 &  + 4,577,802  & + 16.0\\
 	&  	20   & + 479,953  &    + 1.5 &  + 1,809,060  & + 5.8\\
\hline
  	&  	0    & + 231,532  &    + 6.1 &   + 2,103,904  & + 111.0\\
U14 &  	10   & + 391,966  &    + 10.4 &  + 1,184,917  & + 39.8\\
 	&  	20   & + 485,693  &    + 12.9 &  + 494,986   & + 13.2\\
\hline
  	&  	0    & -  76,685  &   - 0.2 &   + 11,622,839  & + 49.5\\
U15 &  	10   & + 524,424  &   + 1.5 &  + 5,459,972   & + 18.1\\
 	&  	20   & + 442,970  &   + 1.3 &  + 3,175,312   & + 9.8\\
\hline
\end{tabular}
\end{center}
\end{table}

\section{Conclusions}

We have presented a new Robust Optimization method for tackling price uncertainty in  energy offering for a price-taker generating company operating in a competitive energy market. Our investigations have been motivated by a review of a Robust Optimization method based on $\Gamma$-Robustness for price-uncertain energy offering proposed in \cite{BaCo11}:
though being an important reference in literature, this method presents a number of issues that may sensibly limit its application in practice and that concretely expose a company to the risk of presenting offering curves resulting into suboptimal and even infeasible accepted offers. To tackle all such issues, we have proposed an alternative $\Gamma$-Robust Optimization method that requires to solve one single robust counterpart, considering an intermediate level of protection between null and full protection, and to make energy offers at zero price, practically eliminating the risk of non-acceptance. Computational results on a set of instances provided by our industrial partners show that our method is able to grant a very high increase in the profit with respect to solutions obtained by a strategy that we have observed among professionals and to solutions obtained by imposing full protection against price deviations.
As a future direction of investigations, we intend to study the adaption of alternative Robust Optimization models that grant a better representation of price uncertainty, such as Multiband Robustness, a refinement and generalization of $\Gamma$-Robustness that exploits histogram-like uncertainty sets (see e.g., \cite{BaEtAl14,BuDA12a}). A better modeling of the uncertain price data would lead indeed to further reduction in the conservatism of robust solutions, and thus increases in the profit, while not reducing protection against price deviations.

\section*{Acknowledgments}
\noindent
This first author was partially supported by the \emph{Einstein Center for Mathematics Berlin} (ECMath) through Project MI4 (ROUAN) and by the \emph{German Federal Ministry of Education and Research} (BMBF) through Project VINO and Project \emph{ROBUKOM} \cite{BaEtAl14}.
The first author was also partially supported by the COST Action TD1207 and
would like to acknowledge networking support by the COST Action TD1207.

\appendices

\section{Reference formulation} \label{sec:appendixRefForm}
\noindent
The set of feasible power solutions $P$ is described by the mixed integer quadratic programming problem introduced below, using the following set of constants and variables characterizing the considered single generation unit.

\noindent
Constants:
\begin{itemize}
  \item $\lambda_t$ market price in period $t \in T$;
  \item $a, b \geq 0$ coefficients of the quadratic and linear term of the quadratic function expressing generation cost;
  \item $c^{F}$ fixed cost of generation in one single period;
  \item $c^{SU}_\tau$ startup cost after $\tau$ consecutive periods of turned-off status;
  \item $P^{\min}, P^{\max}$ minimum and maximum power output;
  \item $R^{\nearrow}, R^{\searrow}$ ramp-up and ramp-down limits;
  \item $R^{SU}, R^{SD}$ startup and shutdown ramp limits;
  \item $U_{i}, D_{i}$ minimum up and down time.
\end{itemize}

\noindent
Decision variables:
\begin{itemize}
  \item $p_{t} \geq 0$ - power output variable for period $t \in T$;
    \item $suc_{t} \geq 0$ - startup cost variable for period $t \in T$;
  \item $u_{t}  \in \{0,1\}$ - status variable for period $t \in T$. A variable $u_{t}$ is equal to 1 if the unit is \emph{on} in period $t$ and 0 if the unit is \emph{off});
  \item $v_{t} \in \{0,1\}$ - startup variable for period $t \in T$. A variable $v_{t}$ is equal to 1 if the unit is \emph{turned on} in period $t$ and 0 otherwise;
  \item $w_{t} \in \{0,1\}$ - shutdown variable for period $t \in T$. A variable $w_{t}$ is equal to 1 if the unit is \emph{turned off} in period $t$ and 0 otherwise.
\end{itemize}

\noindent
We formulate the unit commitment problem as follows:
\begin{align}
\max
&
\sum_{t \in T}
\left\{
\lambda_t p_t
-
    \left[
        a  \left( p_t \right)^2
        + b \hspace{0.1cm}  p_t
        + c^{F}  u_t
        + suc_{t}
    \right]
\right\}
&&
\label{MILP-objFunction}
\\
&
suc_t \geq c^{SU}_\tau \left( u_{t} - \sum_{k = 1}^{\tau} u_{t-k} \right)
\qquad t \in T
\label{MILP-SUCConstraint}
\\
&
P^{\min} \hspace{0.1cm} u_{t} \leq p_{t} \leq P^{\max} \hspace{0.1cm} u_{t}
\qquad t \in T
\label{MILP-variableBoundPowerConstraint}
\\
&
p_{t} \leq p_{t-1} + R^{\nearrow} u_{t}
+ \left(R^{SU} - R^{\nearrow}\right) v_{t}
\qquad
t \in T
\label{MILP-RampUpConstraint}
\\
&
p_{t} \geq p_{t-1} - R^{\searrow} u_{t-1}
+ \left(R^{\searrow} - R^{SD} \right) w_{t}
\hspace{0.4cm}
t \in T
\label{MILP-RampDownConstraint}
\\
&
\sum_{\tau = t - U + 1}^{t} v_{\tau} \leq u_{t}
\qquad
t \in \{U + 1, \ldots, |T|\}
\label{MILP-MinUpTimeConstraint}
\\
&
\sum_{\tau = t - D + 1}^{t} w_{\tau} \leq 1-  u_{t}
\qquad
t \in \{D + 1, \ldots, |T|\}
\label{MILP-MinDownTimeConstraint}
\\
&
w_{t} = v_{t} + u_{t-1} - u_{t}
\qquad
t \in T
\label{MILP-VWConstraint}
\\
&
p_{t} \geq 0
\qquad
t \in T
\label{MILP-PowerVar}
\\
&
suc_{t} \geq 0
\qquad
t \in T
\label{MILP-SUCVar}
\\
&
u_{t} \in \{0,1\}
\qquad
t \in T
\label{MILP-statusVar}
\\
&
v_{t} \in \{0,1\}
\qquad
t \in T
\label{MILP-startupVar}
\\
&
w_{t} \in \{0,1\}
\qquad
t \in T
\label{MILP-shutdownVar}
\end{align}

The objective function aims at maximizing the total profit of the unit over the time horizon, obtained as the sum of the difference of the total revenue and the total cost in each time period. The total cost of generation in one period is equal to the sum of the quadratic cost function plus the startup cost. The constraints \eqref{MILP-SUCConstraint} express the linking between the startup cost variables and the status of the generation unit multiplied by the corresponding startup constants. The constraints \eqref{MILP-variableBoundPowerConstraint} are variable bound constraints connecting the power output variables to the the status variable in each time period. The ramp-up and start-up limits and the ramp limits at startup and shutdown are respectively imposed by the constraints \eqref{MILP-RampUpConstraint} and \eqref{MILP-RampDownConstraint}. The constraints \eqref{MILP-MinUpTimeConstraint} and \eqref{MILP-MinDownTimeConstraint} impose the minimum up and down time of each unit.
Finally, (\ref{MILP-PowerVar}-\ref{MILP-shutdownVar}) declare the decision variables of the problem.

In comparison to the formulation used in \cite{BaCo11}, our reference formulation contains two main modeling and polyhedral refinements:
a) we consider a quadratic cost function $c_t(p_{t})$, which provides a refined representation of the production costs w.r.t. the linear
approximation adopted in \cite{BaCo11}; b) we consider a stronger formulation of the problem, which uses the inequalities providing a complete and compact description of the convex hull of the minimum up and down times of a unit (introduced in \cite{RaTa05}).
We are aware that in literature there exist further ways for strengthening other families of unit commitment constraints (for example, the ramp constraints - see \cite{TaEtAl15} for an overview) and refined linear representation of the quadratic cost function (e.g., by perspective cuts \cite{FraEtAl09l03}). However, since we can solve (\ref{MILP-objFunction}-\ref{MILP-shutdownVar}) practically instantaneously, such refinements resulted not necessary and we decided to not include them.

\ifCLASSOPTIONcaptionsoff
  \newpage
\fi


\begin{thebibliography}{1}

\bibitem{BaCo11}
L.~Baringo, A.~J. Conejo: \emph{Offering Strategy Via Robust Optimization}, IEEE Trans. Power Syst. 26 (3), 418-425, 2011

\bibitem{BaEtAl14}
T. Bauschert, C. B\"using, F. D'Andreagiovanni, A.M.C.A. Koster, M. Kutschka, U. Steglich:: \emph{Network Planning under Demand Uncertainty with Robust Optimization}. IEEE Comm. Mag. 52 (2), 178-185, 2014, DOI: 10.1109/MCOM.2014.6736760

\bibitem{BeElNe09} A.~Ben-Tal, L.~El Ghaoui, A.~Nemirovski: \emph{Robust Optimization}, Heidelberg, Germany: Springer, 2009

\bibitem{BeBrCa11} D.~Bertsimas, D.~Brown, C.~Caramanis: \emph{Theory and Applications of Robust Optimization}, SIAM Review 53 (3), 464-501, 2011

\bibitem{BeSi04} D.~Bertsimas, M.~Sim: \emph{The Price of Robustness}, Oper. Res. 52 (1), 35-53, 2004

\bibitem{BuDA12a} C. B\"using, , F.D'Andreagiovanni: New Results about Multi-band Uncertainty in Robust Optimization. In: Klasing, R. (ed.) Experimental Algorithms - SEA 2012, LNCS, vol. 7276, pp. 63-74. Springer, Heidelberg, 2012, DOI: 10.1007/978-3-642-30850-5\_7

\bibitem{CoEtAl02} A.J. Conejo, F.J. Nogales, J.M. Arroyo: Price-taker bidding strategy under price uncertainty. IEEE Trans. Power Syst. 17 (4), 1081-1088, 2002

\bibitem{CoEtAl04} A.J. Conejo, F.J. Nogales, J.M. Arroyo, R. Garcia-Bertrand: \emph{Risk-Constrained Self-Scheduling of a Thermal Power Producer}. IEEE Trans. Power Syst. 19 (3), 1569-1574, 2004

\bibitem{CPLEX} IBM ILOG CPLEX: http://www-01.ibm.com/software/integration/ optimization/cplex-optimizer

\bibitem{Da93} A. David: \emph{Competitive bidding in electricity supply}. Proc. Inst. Elect. Eng., Gen., Transm., Distrib. C, 140 (5), 421-426, 1993

\bibitem{DaWe00} A.K. Davin, F. Wen: \emph{Strategic bidding in competitive electricity markets: a literature survey}. Proc. IEEE Power Eng. Soc. 4, 2168-2173, 2000

\bibitem{DeEtAl03}
M. Denton, A. Palmer, R. Masiello, P. Skantze: \emph{Managing market risk in energy}. IEEE Trans. Power Syst. 18 (2), 494-502, 2003

\bibitem{FraEtAl09l03}
A. {Frangioni}, C. {Gentile}, F. {Lacalandra}: \emph{Tighter Approximated MILP Formulations for Unit Commitment Problems}. IEEE Trans. Power Syst. 24 (1), 105-113, 2009

\bibitem{GME_website} Gestore Mercati Energetici.
http://www.mercatoelettrico.org/en/


\bibitem{GoBa04}
V. Gountis, A. Bakirtzis: \emph{Bidding strategies for electricity producers
in a competitive electricity marketplace}. IEEE Trans. Power
Syst. 19 (1), 356-365, 2004

\bibitem{HrEtAl13}
H. Pandžic, J.M. Morales, A.J. Conejo, I. Kuzle:
\emph{Offering model for a virtual power plant based on stochastic programming}.
J. Appl. Energy 105, 282-292, 2013

\bibitem{Ja05} R. Jab: \emph{Robust self-scheduling under price uncertainty using conditional value-at-risk}. IEEE Trans. Power Syst. 20 (4), 1852-1858, 2005

\bibitem{KwFr12} R.H. Kwon, D. Frances: \emph{Optimization-Based Bidding in Day-Ahead Electricity Auction Markets: A Review of Models for Power Producers}.  In: Sorokin, A. et al. (eds.) Handbook of Networks in Power Systems, pp. 41-59. Springer, Heidelberg, 2012

\bibitem{LiShLi07}
T. Li, M. Shahidehpour, Z. Li:
\emph{Risk-Constrained Bidding Strategy With Stochastic Unit Commitment}.
IEEE Trans. Power Syst. 22 (1), 449-458, 2007

\bibitem{LiShQu11}
G. Li, J. Shi, X. Qu: \emph{Modeling methods for GenCo bidding strategy optimization in the liberalized electricity spot market – A state-of-the-art review}. Energy
36 (8), 4686-4700, 2011

\bibitem{RaTa05}
Rajan, D., Takriti, S.: \emph{Minimum up/down polytopes of the unit commitment problem with start-up costs}. IBM Research Report RC23628, 2005


\bibitem{TaEtAl15}
M. Tahanan, W. van Ackooij, A. Frangioni, F. Lacalandra:
\emph{Large-scale Unit Commitment under uncertainty}.
4OR 13 (2), 115-171, 2015

\bibitem{YaSh04}
H. Yamin, S. Shahidehpour: \emph{Risk and profit in self-scheduling for
gencos}. IEEE Trans. Power Syst. 19 (4), 2104-2106, 2004

\end{thebibliography}
\end{document}